\documentclass[conference]{IEEEtran}
\usepackage[latin9]{inputenc}
\usepackage{amsmath,amssymb,color}
\usepackage{graphicx,epsfig,framed}
\usepackage{mathptmx,times} 

\makeatletter
\newcommand{\MYfooter}{\smash{
\hfil\parbox[t][\height][t]{\textwidth}{}\hfil\hbox{}}}
\makeatletter
\def\ps@IEEEtitlepagestyle{%
\def\@oddhead{\mbox{}2016 ICSEE International Conference on the Science of Electrical Engineering \rightmark \hfil }

 \def\@oddfoot{\MYfooter{978-1-5090-2152-9/16/\$31.00\:\copyright2016\:IEEE}}
 
\def\@evenfoot{\MYfooter}}
\def\ps@headings{%
\def\@oddhead{\mbox{}2016 ICSEE International Conference on the Science of Electrical Engineering \rightmark \hfil }} \makeatother 
\newtheorem{theorem}{Theorem}[section]

\newtheorem{lemma}[theorem]{Lemma}

\newtheorem{proposition}[theorem]{Proposition}

\newtheorem{remark}[theorem]{Remark}
\numberwithin{equation}{section}

\newcommand{\BE}{\begin{equation}}
\newcommand{\BEQ}[1]{\BE\mathlabel{#1}} 
\newcommand{\EEQ}{\end{equation}}
\newcommand{\rfb}[1]{\mbox{\rm
   (\ref{#1})}\ifx\undefined\stillediting\else:\fbox{$#1$}\fi}

\newcommand{\nline}  {{\mathbb N}}
\newcommand{\rline}  {{\mathbb R}}

\newcommand{\PPP}{{\mathbf P}}

\newcommand{\Gscr} {{\cal G}}
\newcommand{\Kscr} {{\cal K}}

\newcommand{\Mscr} {{\cal M}}

\newcommand{\Sscr} {{\cal S}}
\newcommand{\Tscr} {{\cal T}}

\newcommand{\mm}     {{\hbox{\hskip 0.5pt}}}
\newcommand{\m}      {{\hbox{\hskip 1pt}}}
\newcommand{\nm}     {{\hbox{\hskip -3pt}}}
\newcommand{\bluff}  {{\hbox{\raise 15pt \hbox{\mm}}}}
\newcommand{\sbluff} {{\hbox{\raise  9pt \hbox{\mm}}}}

\renewcommand{\l}    {{\lambda}}

\newcommand{\e}      {{\varepsilon}}

\newcommand{\dd}     {{\rm d\hbox{\hskip 0.5pt}}}

\newcommand{\FORALL} {{\hbox{$\hskip 11mm \forall \;$}}}
\newcommand{\rarrow} {\mathop{\rightarrow}}                  
\newcommand{\sat}    {\mathop{\rm sat}}
\let\oldlabel=\label

\renewcommand{\label}[1]{\leavevmode\smash{\raise 10pt\llap
             {\fbox{\scriptsize#1}}}\oldlabel{#1}}
\newcommand{\mathlabel}[1]{\smash{\raise 9pt\llap
             {\scriptsize(#1)}}\label{#1}}

\renewcommand{\label}[1]{\oldlabel{#1}}
\renewcommand{\mathlabel}[1]{\label{#1}}

\renewcommand{\Re}{{\rm Re\,}}

\newcommand{\bbm}[1]{\left[\begin{matrix} #1 \end{matrix}\right]}

\makeatother

\begin{document}
\title{Stability of the integral control of\\ stable nonlinear
       systems} 

\author{\IEEEauthorblockN{George Weiss} \IEEEauthorblockA{School of
        Electrical Engineering, Tel Aviv University\\ Ramat Aviv,
        Israel, 69978} \and \IEEEauthorblockN{Vivek Natarajan} 
        \IEEEauthorblockA{Systems and Control Engineering Group, 
        Indian Institute\\ of Technology (IIT) Bombay, Mumbai 400076, 
        India}}
\thanks{This work was partially supported by grant no. 800/14
        of the Israel Science Foundation, a grant from the Israeli 
        Ministry of Infrastructure, Energy and Water and a Leverhulme
        Fellowship at the Univ. of Warwick.}
\maketitle

\begin{abstract}
PI controllers are the most widespread type of controllers and there
is an intuitive understanding that if their gains are sufficiently
small and of the correct sign, then they ``always'' work. In this
paper we try to give some rigorous backing to this claim, under
specific assumptions. Let $\PPP$ be a nonlinear system described by
$\dot x=f(x,u)$, $y=g(x)$, where the state trajectory $x$ takes values
in $\rline^n$, $u$ and $y$ are scalar and $f,g$ are of class $C^1$. We
assume that there is a Lipschitz function $\Xi:[u_{min},\m
u_{max}]\rarrow\rline^n$ such that for every constant input
$u_0\in[u_{min},\m u_{max}]$, $\Xi(u_0)$ is an exponentially stable
equilibrium point of \m $\PPP$. We also assume that $G(u)=g(\Xi(u))$,
which is the steady state input-output map of $\PPP$, is strictly
increasing. Denoting $y_{min}=G(u_{min})$ and $y_{max}=G (u_{max})$,
we assume that the reference value $r$ is in $(y_{min},\m
y_{max})$. Our aim is that $y$ should track $r$, i.e., $y\rarrow r$ as
$t\rarrow\infty$, while the input of $\PPP$ is only allowed to be in
$[u_{min},\m u_{max}]$. For this, we introduce a variation of the
integrator, called the saturating integrator, and connect it in
feedback with $\PPP$ in the standard way, with gain $k>0$. We show
that for any small enough $k$, the closed-loop system is (locally)
exponentially stable around an equilibrium point $(\Xi(u_r),u_r)$,
with a ``large'' region of attraction $X_T\subset\rline^n\times
[u_{min},\m u_{max}]$. When the state $(x(t),u(t))$ of the 
closed-loop system converges to $(\Xi(u_r),u_r)$, then the tracking
error $r-y$ tends to zero. The compact set $X_T$ can be made larger by
choosing a larger parameter $T>0$, resulting in smaller $k$.
\end{abstract}

\section{Introduction and the definition of the saturating integrator} 
\label{sec1_bis} 

In this short paper we prove some results about the integral control
of stable nonlinear systems. Let the nonlinear time-invariant system
$\PPP$ be described by \vspace{-1mm}
\BEQ{P_syst}
   \dot x \m=\m f(x,u) \m,\qquad y \m=\m g(x) \m,
\end{equation}
where $f$ and $g$ are $C^1$ functions. The state of this system is
$x\in\rline^n$, the input $u$ and the output $y$ are scalar. We assume
that for each constant input function $u_0$ in a certain range
$[u_{min},\m u_{max}]$, $\PPP$ has a locally asymptotically stable
equilibrium point $\Xi(u_0)$ and the function $\Xi:[u_{min}, \m
u_{max}]\rarrow\rline^n$ is Lipschitz continuous. We are not allowed
to apply to $\PPP$ an input function with values outside the range
$[u_{min},\m u_{max}]$, either because the system may become unstable,
or because of actuator saturation, or because of safety considerations
(such as overvoltage on components) - the reason for this limitation
is not relevant for the theory developed here. More technical
assumptions will be stated in the later sections, here we want to
explain the idea.

It is intuitively appealing to regard $\PPP$ as being approximately
modelled by the memoryless system $y=g(\Xi(u))$, and this would be 
close to correct if $u$ were a very slowly changing signal with 
values in the range $[u_{min},\m u_{max}]$. We assume that the
function $G=g\circ\Xi$ is strictly increasing and we denote
$$ y_{min} \m=\m G(u_{min}) \m,\qquad y_{max} \m=\m
   G(u_{max}) \m.$$

The control objective is to make $y$ track a constant (but not given
a-priori) reference signal $r\in(y_{min},\m y_{max})$, while not
allowing the input signal to exit the range $[u_{min},\m u_{max}]$. If
$\PPP$ is replaced with the memoryless model $y=G(u)$ mentioned above,
then this control objective can be achieved using an integral
controller with saturation: \vspace{-1mm}
$$\dot v(t) \m=\m k [r-y(t)] \m,\qquad u(t) \m=\m \sat(v(t)) \m,$$
where $\sat$ denotes a saturation function that does not allow $u$ to
exit the range $[u_{min},\m u_{max}]$, and $u=v$ if $v$ is inside the
allowed range. It is not difficult to show that (for the memoryless
model) this would work, i.e., the closed-loop system would be stable
and we would have $y(t)\rarrow r$ as $t\rarrow\infty$.

The above very simple result (for the memoryless model) can be shown
using a quadratic Lyapunov function, or it may be regarded as an
application of the famous circle criterion, for which we refer to the
nice survey \cite{JaLoRy:11}. Even for this situation, the saturation
as described is not satisfactory, because during a fault the state $u$
of the integrator may reach a very large value (a phenomenon called
``windup''), from which it would take a long time to recover after the
fault. A better way to build the integrator is to prevent its state
from exiting the range $[u_{min},\m u_{max}]$. There are different
ways to do this, and such controllers are said to have {\em
anti-windup}. There is a rich literature on control with anti-windup,
with a much wider meaning for the concept, see for instance 
\cite{KoCaMoNe:94}, \cite{WuLu:04}, \cite{ZaccTeel:02}. We propose one
very particular controller with anti-windup, which we call the {\em 
saturating integrator}, a dynamical system defined by \vspace{-1mm}
\BEQ{sat_int_bis}
   \dot u \m=\m \Sscr(u,w) \m,\vspace{-2mm}
\end{equation}
where \vspace{-2mm}
\BEQ{Sadiq_Khan_bis}
   \Sscr(u,w) \m=\m \begin{cases}
   \m w^+ \ \ &\mbox{ if }\ \ u\leq u_{min}\m,\\
   \m w   \ \ &\mbox{ if }\ \ u\in(u_{min},u_{max}) \m,\\ 
   \m w^- \ \ &\mbox{ if }\ \ u\geq u_{max} \m.
\end{cases} \end{equation}
Here $w^+$ is the positive part of $w$ and $w^-$ is the negative
part of $w$: \vspace{-2mm}
$$w^+ \m=\m \max\{w,0\} \m,\qquad w^- \m=\m \min\{w,0\} \m.$$
The state of the saturating integrator is $u$ and its state 
space is the interval $[u_{min},\m u_{max}]$. 

If $w$ is a continuous function with finitely many zeros in every
finite interval, then it is easy to define the corresponding state
trajectories of the saturating integrator, even though the function
$\Sscr$ is not continuous. However, if the zeros of $w$ have an
accumulation point, then the definition of state trajectories $u$ of
this system becomes problematic. For instance, if $u(0)=u_{max}$ and
$w(t)=t \sin(1/t)$, then it is not obvious what the function $u$ is.
To overcome this problem, let us first consider only inputs $w$ that
are not problematic, for instance, polynomials. It is easy to check
that if $u_1$ and $u_2$ are state trajectories of the saturating 
integrator corresponding to the polynomials $w_1$ and $w_2$,
respectively (and any initial states), and at some moment $t\geq 0$
we have $u_2(t)\geq u_1(t)$, then \vspace{-1mm}
$$\frac{\dd}{\dd t}[u_2(t)-u_1(t)] \m\leq\m w_2(t)-w_1(t) \m,$$
which implies that \vspace{-1mm}
$$\frac{\dd}{\dd t}|u_2(t)-u_1(t)| \m\leq\m |w_2(t)-w_1(t)| \m,$$
and by a symmetric argument this last inequality is true also when
$u_2(t)<u_1(t)$. It follows that \vspace{-1mm}
$$ |u_2(t)-u_1(t)| \m\leq\m |u_2(0)-u_1(0)|+\int_0^t|w_2(\sigma)-
   w_1(\sigma)| \dd\sigma \m.$$
This shows that $u(t)$ from \rfb{sat_int_bis} depends Lipschitz
continuously both on $u(0)$ and also on $w$ considered with the $L^1$
norm. Indeed, for $u_2(0)=u_1(0)$ we can write the last estimate as
\BEQ{Jim_and_Helen_barbeque}
   |u_2(t)-u_1(t)| \m\leq\m \|w_2-w_1\|_{L^1[0,t]} \m.
\end{equation}
Hence, by continuous extension, we can define $u(t)$ for any
input $w\in L^1[0,t]$ (because the polynomials are dense in
$L^1[0,t]$).  In block diagrams (such as Figure 1) we use the symbol
$\int\nm\Sscr$ to denote the saturating integrator. The saturating
integrator has been used also in \cite{LiWe:15}.

The main results of this paper concern the feedback system shown
in Figure 1, which is described by \rfb{P_syst}, \rfb{sat_int_bis}
and $w=k(r-y)$. The state of the closed-loop system is $(x(t),u(t))$
and its state space is \vspace{-4mm}
\BEQ{X}
   \m\ \ \ X \m=\m \rline^n \times [u_{min},\m u_{max}] \m.
\end{equation}

\m\vspace{-5mm}
$$\includegraphics[scale=0.15]{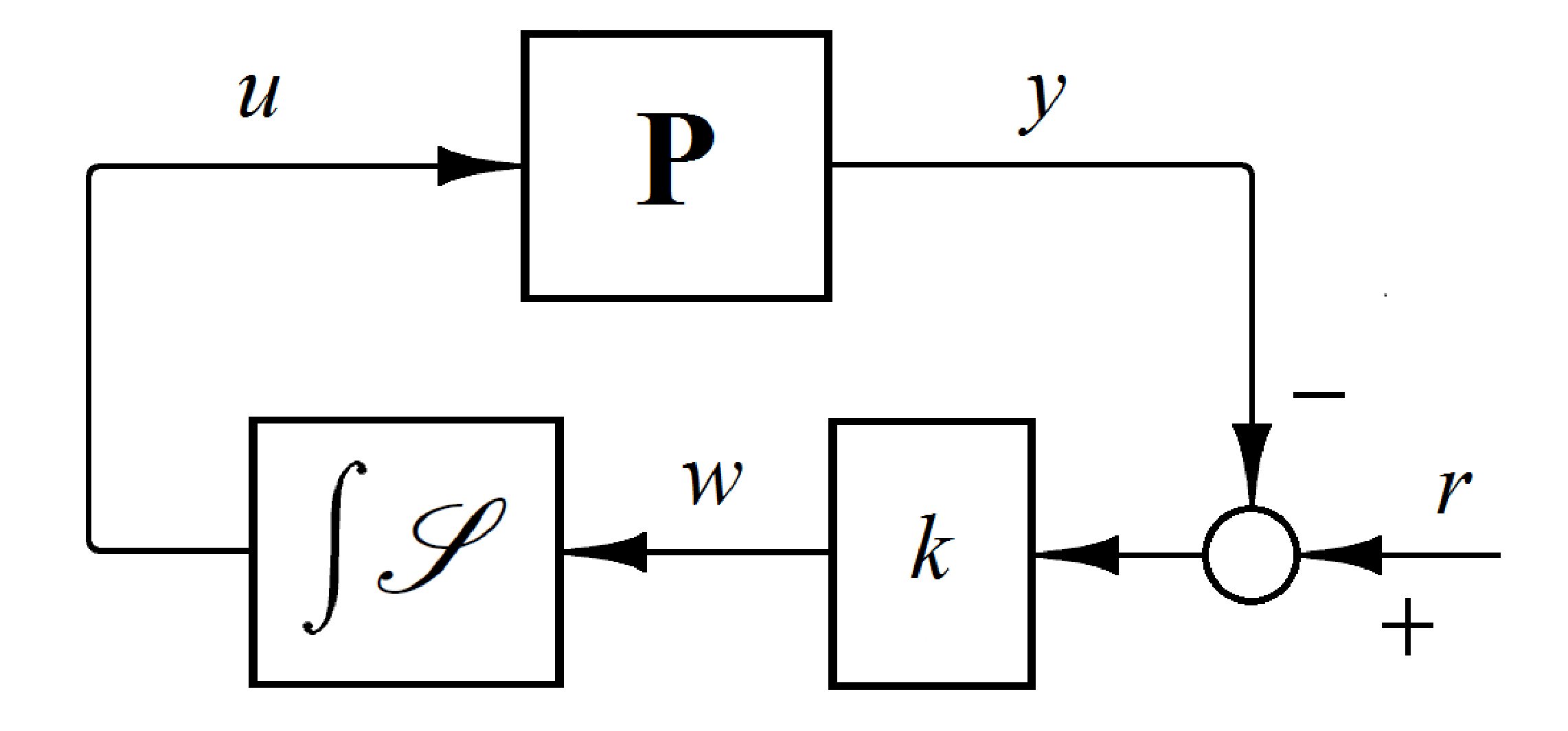}$$
\centerline{ \parbox{3.4in}{Figure 1. The closed-loop system formed
   from the plant $\PPP$, the saturating integrator $\int\nm\Sscr$
   and the constant gain $k>0$, with the constant reference $r$.}}
\medskip

An informal statement of the main result of this paper is that with
the saturating integrator as the controller in the feedback loop, 
under reasonable assumptions on the plant, for any constant 
reference $r$ in the range $(y_{min},\m y_{max})$, the following 
holds: For any small enough feedback gain $k>0$, the closed-loop
system shown in Figure 1 is locally asymptotically stable around an 
equilibrium point, with a ``large'' region of attraction. When the
state converges to this equilibrium point, then the tracking error
$r-y$ tends to zero. 

The precise statement of the main results and their proof will be
given in Sections \ref{sec3_bis} and \ref{sec4_bis}. This theory has
been developed with a very specific example in mind: the control of
the virtual field current in a synchronverter. Explaining the context
of that application would take several pages and instead we just refer
to the papers \cite{NaWe:15,NaWe:16}. The material in this paper was
originally meant to be a lemma in \cite{NaWe:16}, but then it grew too
long. The authors believe that the results are relevant for many more
applications and they are amenable to various generalizations.

Our main results are related to those in \cite{FlLoRy:03}, where
$\PPP$ is assumed to be built from a stable linear system connected in
cascade with nondecreasing nonlinear functions (memoryless systems)
both at its input and at its output. It seems that for such $\PPP$,
our Theorem \ref{Tamir} follows from Theorem 7 and Remark 8 in
\cite{FlLoRy:03}. Another class of related results concerns the
situation when $\PPP$ is assumed to be impedance passive, which allows
an entirely different approach to the proof of set point regulation 
with closed-loop stability, with arbitrary positive gain, see for
instance \cite{Orte:1998}, \cite{OrRo:12}, \cite{JaWe:09}.

\section{Nonlinear systems with slowly varying inputs} 
\label{sec2_bis} 

In this section we investigate the behaviour of a nonlinear system
$\PPP$ from \rfb{P_syst} (with $f$ and $g$ of class $C^1$), assuming
that it has certain stability properties formulated in Assumption 1
below. It is well-known that for any initial state $x(0)$ and any
continuous input function $u$, the differential equation in
\rfb{P_syst} has a unique solution defined on some maximal interval
$[0,t^*)$ (possibly $t^*=\infty$), see \cite[Appendix C] {Sont:1990}
or \cite[Chapter 3]{Kha:02} for good discussions of this topic. It is
important to note that if $t^*$ is finite, then $\limsup_{t\rarrow
t^*}\|x(t)\|=\infty$, see for instance Exercise 3.26 in \cite{Kha:02}
(see also Corollary 2.3 in \cite{JaWe:09}). In this case, we say that
the state trajectory has a {\em finite escape time} $t^*$. The two
lemmas in this section imply that certain state trajectories of $\PPP$
remain bounded as long as they exist, and this of course implies that
they exist for all $t\geq 0$. In many places, our arguments should
contain phrases like ``if the solution exists for this $t$, then
...''. However, in order to make this text less clumsy, we will
discuss about these state trajectories as if it is clear from the
start that they exist for all $t\geq 0$.

{\bf Notation.} \m For any interval $J$, any $\alpha>0$ and any
$m\in\nline$, we denote by Lip$_\alpha(J;\rline^m)$ the set of those
$u:J\rarrow\rline^m$ which are Lipschitz continuous with Lipschitz
constant $\alpha$. If $u$ is defined on a larger set containing $J$,
then $u\in$\m\m Lip$_\alpha(J;\rline^m)$ means that the restriction of
$u$ to $J$ is in Lip$_\alpha(J;\rline^m)$. \vspace{-2mm}

\begin{framed} \vspace{-2mm}
{\bf Assumption 1.} There exist real numbers $u_{min}<u_{max}$, 
$\alpha>0$ and a function $\Xi\in$ Lip$_\alpha(([u_{min},u_{max}];
\rline^n)$ such that  
$$f(\Xi(u),u) \m=\m 0 \FORALL u\in [u_{min},u_{max}] \m,$$
i.e., for each $u_0\in[u_{min},u_{max}]$, $\Xi(u_0)$ is an 
equilibrium point that corresponds to the constant input $u_0$.

Moreover, $\PPP$ is {\em uniformly exponentially stable} around these
equilibrium points. This means that there exist \hbox{$\e_0>0$},
$\l>0$ and $m\geq 1$ such that for each constant input function
$u_0\in [u_{min},\m u_{max}]$, the following holds:

If \ $\|x(0)-\Xi(u_0)\|\leq\e_0$, then for every $t\geq 0$, 
\BEQ{exp_decay}
   \|x(t)-\Xi(u_0)\| \m\leq\m m e^{-\l t}\|x(0)-\Xi(u_0)\| \m.
\vspace{-3mm} \end{equation}
\end{framed} \vspace{-2mm}

\begin{remark}
The uniform exponential stability condition above can be checked by
linearization: If the Jacobian matrices \vspace{-2mm}
$$ A(u_0) \m=\m \left. \frac{\partial f(x,u)}{\partial x} \right|_
   {\nm\begin{array}{c} \scriptstyle x=\Xi(u_0)\vspace{-1mm}\\ 
   \scriptstyle u=u_0 \end{array}} \m\in\m \rline^{n\times n}
   \vspace{-1mm}$$
have eigenvalues bounded away from the right half-plane, \vspace{-1mm}
$$ \max \Re \sigma(A(u_0)) \m\leq\m \l_0 \m<\m 0 \FORALL u_0\in
   [u_{min},\m u_{max}] \m,\vspace{-1mm}$$
then $\PPP$ is uniformly exponentially stable, see (11.16) in 
\cite{Kha:02}. Under Assumption 1, \m $\max\Re\sigma(A(u_0))$ is a 
continuous function of $u_0$. Hence, if this function is always
negative, then by the compactness of $[u_{min},\m u_{max}]$, its
maximum is also negative. Thus, for the uniform exponential stability
we only have to check that each of the matrices $A(u_0)$ is stable. 
\end{remark} \smallskip

The following two lemmas show that under the above assumption, if the
input $u$ changes sufficiently slowly and stays in the relevant range
of values $[u_{min},u_{max}]$, and if $x(0)$ is close to its momentary
equilibium value $\Xi(u(0))$, then for all times \m $t>0$, $x(t)$
remains close to $\Xi(u(t))$. These results are related to those in
Section 9.3 of \cite{Kha:02}, where the proof technique is to
construct and use Lyapunov functions. It is difficult to see the 
precise relationship between our lemmas and those in \cite{Kha:02}.

\medskip {\color{blue}
\begin{lemma} \label{slow_input}
Assume that $\PPP$ satisfies Assumption 1. 
 
Then there exists $\kappa>0$ and $T>0$ such that for every $\e\in[0,
\e_0]$ and every \m $u\in{\rm Lip}_{\kappa\e}([0,\infty);\rline)$ 
with \m $u(t)\in[u_{min},u_{max}]$ \m for all $t\geq 0$, the 
following holds: for all $t\geq T$, \vspace{-2mm}
\BEQ{Jutta_Strake}
   \|x(0)-\Xi(u(0))\| \m\leq\m \e \ \Longrightarrow \ 
   \|x(t)-\Xi(u(t))\| \m\leq\m \frac{2}{3}\e \m.
\end{equation}
\end{lemma}}

The proof is in the journal version of this work \cite{WeNa:16}.

\medskip {\color{blue}
\begin{lemma} \label{RenanKhalifa}
Suppose that Assumption 1 holds and let $\kappa,T$ be the positive
constants whose existence was proved in Lemma \ref{slow_input}. Let
$\e\in(0,\e_0]$ and assume that the initial state $x(0)$ and and the
input $u$ of $\PPP$ satisfy
$$ \|x(0)-\Xi(u(0))\| \m\leq\m \e \m,\qquad u\in{\rm Lip}_{\kappa
   \e}([0,\infty);\rline)$$
and \m $u(t)\in[u_{min},\m u_{max}]$ \m for all $t\geq 0$. Then
\vspace{-2mm}
\BEQ{Ben_plays_steam}
   \|x(t)-\Xi(u(t))\| \m<\m \left( m+\frac{1}{6} \right) \e 
   \FORALL t\geq 0 \m.
\end{equation}
\end{lemma}}

For the proof we refer again to the journal version \cite{WeNa:16}.

\section{The closed-loop system} 
\label{sec3_bis} 

In this section we discuss the behaviour of the closed-loop system
from Figure 1, with the state space $X$ from \rfb{X}. There is an
interesting and somewhat unclear connection between the lemmas in this
section and Tikhonov's theorem concerning singularly perturbed systems
of differential equations, see Theorem 11.1 (and also Theorem 11.2) in
\cite{Kha:02}.  Our system may be regarded as a variation of a
subclass of the systems studied in the cited theorems. However,
Tikhonov's theorem concerns the asymptotic behaviour of the solutions
on a parameter, which in our case is $k$, as $k\rarrow 0$, and this is
different from our concerns (we want to establish stability for a
fixed $k$).

We start with a proposition about local existence and uniqueness of
state trajectories, which is not an obvious fact due to the
discontinuity of $\Sscr$. Note that in this proposition we do not
impose Assumption 1 on $\PPP$.

\medskip {\color{blue}
\begin{proposition} \label{BP_shares_missing}
Let $\PPP$ be described by \rfb{P_syst} with $f$ and $g$ of class 
$C^1$ and let $\int\nm\Sscr$ be the saturating integrator as in 
\rfb{sat_int_bis} and \rfb{Sadiq_Khan_bis}. For every $x_0\in
\rline^n$, every $u_0\in[u_{min},\m u_{max}]$, every $k\geq 0$ 
and every $r\in\rline$ there exists a $\tau>0$ such that the 
closed-loop system from Figure 1 has a unique state trajectory 
$(x,u)$ defined on $[0,\tau)$, such that $x(0)=x_0$ and $u(0)=u_0$.

If $\tau$ is maximal (i.e., the state trajectory cannot be continued
beyond $\tau$) then $\limsup_{t\rarrow\tau}\|x(t)\|=\infty$.
\end{proposition}}

\medskip
For the proof we refer to the journal version of this work.

\begin{framed} \vspace{-2mm}
{\bf Assumption 2.} The system $\PPP$ satisfies Assumption 1 and
moreover, the function
$$ G(u) \m=\m g(\Xi (u)) \m, \qquad u\in [u_{min},\m u_{max}]$$
satisfies the following: There exists $\mu>0$ such that for any
$u_1,u_2\in [u_{min},\m u_{max}]$ with $u_1>u_2$,
\BEQ{Aleppo}
   G(u_1)-G(u_2) \m\geq\m 2\mu(u_1-u_2) \m.
\end{equation}
(If $G$ is differentiable then this is equivalent to $G'\geq 2\mu$.)
\vspace{-2mm} \end{framed}

{\bf Notation.} \m Recall from Section \ref{sec1_bis} that we denote
$y_{min}=G(u_{min})$ and $y_{max}=G(u_{max})$, so that clearly 
$y_{min}<y_{max}$. For any $r\in(y_{min},\m y_{max})$ we define $u_r
=G^{-1}(r)$ and we define $\Gscr_r:[u_{min}-u_r,u_{max}-u_r]\rarrow
\rline$ by shifting the graph of $G$: \vspace{-1mm}
$$\Gscr_r(v) \m=\m G(v+u_r)-r \m,$$
so that $\Gscr_r$ is an increasing Lipschitz function and 
$\Gscr_r(0)=0$. It is clear that \rfb{Aleppo} holds with $\Gscr_r$ in
place of $G$.

\smallskip {\color{blue}
\begin{lemma} \label{KenLivinstone}
Consider the closed-loop system from Figure 1, where $\PPP$ satisfies
Assumption 2, $k>0$ and $r\in(y_{min},\m y_{max})$. Assume that $u(0)
\in[u_{min},\m u_{max}]$ and let $x(0)\in\rline^n$ and $\tau,\eta^*>0$
be such that the closed-loop state trajectory $(x,u)$ exists for $t\in
[0,\tau]$ (and possibly also later) and \vspace{-1mm}
\BEQ{tickets}
   |y(t)-G(u(t))| \m\leq\m \eta^* \FORALL t\in[0,\tau] \m.
\end{equation}

Then for all $t\in[0,\tau]$ we have
\BEQ{Luton_parking}
   |G(u(t))-r| \m\leq\m \max\left\{ \left|\Gscr_r (e^{-\mu k\mm t}
   (u(0)-u_r))\right|,\ 2\eta^* \right\} \m.
\end{equation}
\end{lemma}}

For the proof we refer to the journal version \cite{WeNa:16}.

\smallskip {\color{blue}
\begin{lemma} \label{Donald_Trump}
Consider the closed-loop system from Figure 1, where $\PPP$ satisfies 
Assumption 2 and $r\in(y_{min},\m y_{max})$. Recall the constant 
$\kappa>0$ from Lemma \ref{slow_input}. We introduce the tubular open
set \vspace{-1mm}
$$ W \m=\m \left\{ \xi\in\rline^n \ \left|\ \min_{u_0\in[u_{min},
   u_{max}]} \|\xi-\Xi(u_0)\| < \left(m+\frac{1}{6} \right) \e_0
   \right. \right\} \m.$$
Let $\delta_g$ be a Lipschitz bound of $g$ over $W$. Choose
$\tilde\l,k>0$ such that \vspace{-1mm}
\BEQ{TedCruz}
   \tilde\l \m=\m 2\delta_g \left( m+\frac{1}{6} \right) \m,\qquad 
   k \m<\m \frac{2\kappa}{\delta_g (6m+1)} \m.
\end{equation}
Then there exists $\tau>0$ with the following property: If \m $\e\in
[0,\e_0]$, $u(0)\in [u_{min},\m u_{max}]$ and
\BEQ{Iona}
   \|x(0)-\Xi(u(0))\| \m\leq\m \e \m,\qquad 
   |G(u(0))-r| \m\leq\m \tilde\l\e \m,
\end{equation}
then the state trajectory of the closed-loop system exists for all 
$t\geq 0$ and for all $t\geq\tau$ we have
$$ \|x(t)-\Xi(u(t))\| \m\leq\m \frac{2}{3} \e \m,\qquad 
   |G(u(t))-r| \m\leq\m \frac{2}{3} \tilde\l\e \m.$$
\end{lemma}}

For the proof we refer to the journal version of this work.

\medskip {\color{blue}
\begin{theorem} \label{Tamir}
We work under the assumptions of Lemma \ref{Donald_Trump} up to 
and including \rfb{TedCruz}. Then $(\Xi(u_r),u_r)$ is a locally
asymptotically stable equilibrium point of the closed-loop system
from Figure 1, with the state space $X$ from \rfb{X}.

If the initial state $(x(0),u(0))\in X$ of the closed-loop system 
satisfies $\|x(0)-\Xi(u(0))\|\leq\e_0$, then \vspace{-1mm}
$$ x(t) \rarrow \Xi(u_r) \m,\qquad u(t) \rarrow u_r \m,\qquad
   y(t) \rarrow r \m, \vspace{-1mm}$$
and this convergence is at an exponential rate.
\end{theorem}}

\medskip
{\em Proof.} \m We introduce the coordinate transformation
\vspace{-1mm}
$$ \Tscr : X \rarrow \rline^n \times [y_{min}-r,\m y_{max}-r]
   \vspace{-2mm}$$
as follows: \vspace{-1mm}
$$ \bbm{\xi\\ w} \m=\m \Tscr \left( \bbm{x\\ u} \right) \m=\m
   \bbm{x-\Xi(u)\\ G(u)-r} \m.$$
This transformation is invertible, its inverse is
$$ \bbm{x\\ u} \m=\m \Tscr^{-1} \left( \bbm{\xi\\ w} \right) \m=\m
   \bbm{\xi+\Xi \left( \Gscr_r^{-1}(w)+u_r \right) \\ 
   \Gscr_r^{-1}(w)+u_r} \m.$$
Both $\Tscr$ and $\Tscr^{-1}$ are Lipschitz. Note that in the new 
coordinates, the equilibrium point under discussion is $(0,0)$.

Lemma \ref{Donald_Trump} says that there exists a $\tau>0$ such that
for any $\e\in[0,\e_0]$, if the initial state (in the new 
coordinates) is in the rectangular box $\|\xi(0)\|\leq\e$, $|w(0)|
\leq\tilde\l\e$, then for all $t\geq\tau$ the state $(\xi(t),w(t))$
will be in a rectangular box that is $2/3$ times smaller. Clearly 
this implies that (in the new coordinates) the origin is a locally 
asymptotically stable equilibrium point. Moreover, the state 
converges to this equilibrium point at an exponential rate. Clearly
the same conclusions hold for the equilibrium point $(\Xi(u_r),u_r)$
in the original coordinates.

Finally, suppose that the initial state satisfies $u(0)\in[u_{min},\m
u_{max}]$ and $\|x(0)-\Xi(u(0))\|\leq\e_0$. The Lipschitz bound of
$g$, denoted $\delta_g$, can be chosen as large as needed, so that
$\tilde\l$ (given by \rfb{TedCruz}) becomes sufficiently large so that
$|G(u(0))-r|\leq\tilde\l\e_0$ holds. Then we can apply our earlier
argument to conclude that $(x(t),u(t))$ converges to $(\Xi(u_r),u_r)$
at an exponential rate. Since $y(t)=g(x(t))$ and $g$ is a $C^1$
function, it follows that $y(t)$ converges to $g(\Xi(u_r))=G(u_r)=r$
at an exponential rate. \hfill $\Box$ \m

\section{Finding a large domain of attraction} 
\label{sec4_bis} 

In this section we show that, under a well-posedness assumption for
the closed-loop system from Figure 1, we can find a large domain of
attraction for the asymptotically stable equilibrium point whose
existence was proved in Theorem \ref{Tamir}. The following assumption
is stronger than the local well-posedness result in Proposition
\ref{BP_shares_missing}.

\begin{framed} \vspace{-2mm}
{\bf Assumption 3.} There exists $k_0>0$ such that for any $k\in 
[0,k_0]$, the closed-loop system formed by $\PPP$ and the saturating 
integrator, as shown in Figure 1, with any $r\in(y_{min},\m 
y_{max})$, has a unique state trajectory in forward time on the
interval $[0,\infty)$, for any initial state in $X$.

Moreover, at any time $t\geq 0$, the state $(x(t),u(t))$ depends
continuously on the initial state $(x(0),u(0))$.
\vspace{-2mm} \end{framed}

The above assumption is not trivial, because the differential
equations describing the closed-loop system are not continuous (the
discontinuity is in $\Sscr$). It is worth noting that the saturating
integrator is irreversible (in time) and hence the closed-loop system
usually has no uniquely defined backwards (in time) state
trajectories.

\medskip {\color{blue}
\begin{theorem} \label{birthday_of_Kay}
Assume that $\PPP$ satisfies Assumption 2 and moreover, it has
well-defined backwards state trajectories for all $t<0$, corresponding
to any initial state in $\rline^n$ and any constant input in
$[u_{min},\m u_{max}]$. Further, assume that the closed-loop system
from Figure 1 satisfies Assumption 3, and $r\in(y_{min},\m
y_{max})$. Let $T>0$ and define the set $X_T\subset X$ as follows:
$(x_0,u_0)\in X$ belongs to $X_T$ if the state trajectory $z$ of
$\PPP$ starting from $z(0)=x_0$, with constant input $u_0$, satisfies
$\|z(T)-\Xi(u_0)\|\leq\e_0/2$.

Then there exists $k_T\in(0,k_0]$ such that for any $k\in(0,k_T]$, if
the initial state of the closed-loop system is in $X_T$, then the
state trajectory $(x,u)$ of the closed-loop system satisfies
\vspace{-1mm}
$$ x(t) \rarrow \Xi(u_r) \m,\qquad u(t) \rarrow u_r \m,\qquad
   y(t) \rarrow r \m, \vspace{-1mm}$$
and this convergence is at an exponential rate.
\end{theorem}}

\medskip
{\em Proof.} \m Let $\tilde X_T\subset X$ consist of all the points in
the state space $X$ that a state trajectory of the closed-loop system
can reach at some time $t\in[0,T]$, using any fixed value of $k\in 
[0,k_0]$ and starting from an initial state $(x_0,u_0)\in X_T$ at 
time $0$. Obviously $X_T\subset\tilde X_T$. We claim that $\tilde X_T$
is compact.

To prove this claim, first we note that $X_T$ is compact. Indeed, the
system $\PPP$ together with a generator of constant inputs, together
described by the differential equations
$$ \dot x(t) \m=\m f(x(t),u(t)) \nm\qquad \dot u(t) \m=\m 0 \m,$$
has well defined backward state trajectories, given by a continuous
backward flow. The set $X_T$ is the image of \vspace{-1mm}
$$ \Mscr \m=\m \left\{ (z_0,u_0) \in X \ |\ \|z_0-\Xi(u_0)\|\leq
   \frac{\e_0}{2} \right\}$$
through the backward flow mentioned earlier, at time $-T$. It is 
easy to see that $\Mscr$ is compact, and hence $X_T$ (its image
through the backward flow) is also compact.

Now consider the system with state $(x,u,k)$ and state space
$X\times [0,k_0]$ defined by the equations \vspace{-2mm}
$$ \dot x(t) \m=\m f(x(t),u(t)) \m,\quad y(t) \m=\m g(x(t),u(t)) 
   \m,\vspace{-1mm}$$
$$ \dot u(t) \m=\m \Sscr(u(t),k(t)[r-y(t)]) \m,\quad
   \dot k(t) \m=\m 0 \m,\vspace{-1mm}$$
$$ x(0) \m=\m x_0 \m,\ \ \ u(0) \m=\m u_0 \m,\ \ \ 
   k(0) \m=\m k_0 \m.\vspace{-1mm}$$
In other words, this is just the usual closed-loop system, but we
regard $k$ as a constant state variable, that may also take the 
value 0 (which corresponds to constant $u$). From Assumption 3 we
see that this system has a continuous semiflow 
$$ \Phi : X\times[0,k_0]\times[0,\infty)\rarrow X\times[0,k_0]
   \m,\vspace{-1mm}$$
so that $\Phi(x_0,u_0,k_0,t)$ is its state at time $t$. Notice that
\vspace{-1mm}
$$ \tilde X_T \m=\m \Pi\m \Phi(X_T \times [0,k_0] \times [0,T])
   \m,\vspace{-1mm}$$
where $\Pi$ denotes projection onto the first component in the 
product $X\times[0,k_0]$. This implies that indeed $\tilde X_T$ is
compact.

Take $(x_0,u_0)\in X_T$ and let $z$ be the state trajectory of $\PPP$
starting from $z(0)=x_0$, with constant input $u_0$, so that by
assumption $\|z(T)-\Xi(u_0)\|\leq\e_0/2$. Let $(x,u)$ be the state
trajectory of the closed-loop system with some $k\in(0,k_0]$ (to be
specified later) starting from $(x_0,u_0)$. By definition, we know 
that $x(t)\in\tilde X_T$ for all $t\in[0,T]$. We have $u\in{\rm Lip}
_\delta([0,T];\rline)$ where (using the definition of the saturating 
integrator) the Lipschitz bound $\delta$ can be estimated as
\vspace{-1mm} 
$$ \delta \m=\m \max\{k|r-y(t)|\ |\ t\in[0,T]\}$$
$$ \m\leq\m k\max\{|r-g(\xi)|\ |\ (\xi,w)\in\tilde X_T\} \m.$$
Using an argument from the proof of Lemma \ref{slow_input} (see 
\cite{WeNa:16}), we can show that \vspace{-2mm}
$$ \|x(T)-z(T)\| \m\leq\m \frac{L_2\delta T}{L_1} \left[ e^{L_1 T}-1
   \right] \m,$$
where $L_1$ is the Lipschitz constant of $f$ with respect to its first
argument $\xi$, and $L_2$ is the Lipschitz constant of $f$ with
respect to its second argument $w$, when $(\xi,w)\in\tilde X_T$.
Combining the last two estimates, we see that there exists a $p(T)>0$
independent of the initial state in $X_T$ such that $\|x(T)-z(T)\|\leq
p(T)\cdot k$. Thus, we can choose $k_T^1\in(0,k_0]$ small enough so
that $\|x(T)-z(T)\|\leq\e_0/4$ for all $k\in[0,k_T^1]$. Combining this
with $\|z(T)-\Xi(u_0)\|\leq\e_0/2$, we obtain that \vspace{-1mm} 
$$\|x(T)-\Xi(u_0)\| \m\leq\m \frac{3\e_0}{4}\FORALL k\in[0,k_T^1]\m.$$

Finally, it is clear that $|u_0-u(T)|\leq\delta T$, whence
(remembering the constant $\alpha$ from Assumption 1) $\|\Xi(u_0)-\Xi
(u(T))\|\leq\alpha\delta T$. Hence, we can find $k_T^2\in(0,k_T^1]$ 
such that for $k\in[0,k_T^2]$ we have $\|\Xi(u_0)-\Xi(u(T))\|\leq\e_0/
4$. Combining this with the previous estimate, we obtain that
\vspace{-1mm}
$$\|x(T)-\Xi(u(T))\| \m\leq\m \e_0 \FORALL k\in[0,k_T^2] \m.$$
Now we can apply Theorem \ref{Tamir} (starting with the initial time 
$T$) to conclude that for any gain $k\in(0,k_T^2]$ which in addition
satisfies \rfb{TedCruz}, the functions $x,u$ and $y$ converge as
stated. \hfill $\Box$ \m

\begin{remark} 
The reason why we may call $X_T$ a ``large'' domain of attraction is
the following: If $\PPP$ happens to be globally asymptotically stable
(GAS) for every constant input $u_0\in[u_{min},\m u_{max}]$, then {\em
every initial state} of the closed-loop system is contained in one of
the sets $X_T$, if we choose $T$ large enough. If we choose a ``region
of interest'' $\Kscr\subset X$ that is compact, then there exists a 
$k>0$ such that all the closed-loop state trajectories starting from 
$\Kscr$ will converge to the unique equilibrium point. Indeed, the 
interiors of the sets $X_T$ are an open covering of $\Kscr$, so that 
$\Kscr\subset X_T$ if $T$ is large enough. Then we have to choose a
gain $k\leq k_T$. Of course, the price for choosing a very large $T$
is that we may have to choose a very small gain $k$, and this may 
deteriorate the dynamic response of the closed-loop system. 
\end{remark}



\begin{thebibliography}{99}{

\bibitem{FlLoRy:03}
 T.~Fliegner, H.~Logemann and E.P.~Ryan, \m Low-gain integral
 control of continuous-time linear systems subject to input and
 output nonlinearities, {\em Automatica}, vol.~39, 2003,
 pp.~455-462. 

\bibitem{JaLoRy:11}
 B.~Jayawardhana, H.~Logemann and E.P.~Ryan, \m The circle
 criterion and input-to-state stability, {\em IEEE Control Systems
 Magazine}, vol.~31, 2011, pp.~32-67.

\bibitem{JaWe:09}
 B. Jayawardhana and G. Weiss, \m State convergence of passive
 nonlinear systems with an $L^2$ input, {\em IEEE Trans. Automatic
 Control}, vol.~54, 2009, pp.~1723-1727.

\bibitem{Kha:02}
 H.K.~Khalil, \m {\em Nonlinear Systems} (third edition), \m
 Prentice Hall, New Jersey, 2002.

\bibitem{KoCaMoNe:94}
 M.V.~Kothare, P.J.~Campo, M.~Morari, and C.N.~Nett, \m A unified
 framework for the study of anti-windup designs, {\em Automatica},
 vol. 30, 1994, pp.~1869-1883.

\bibitem{KZRK:2015}
 G.~Konstantopoulos, Q.-C.~Zhong, B.~Ren and M.~Krstic, \m Bounded
 integral control for regulating input-to-state stable non-linear
 systems, {\em Proc. of the American Control Conference} (ACC 2015),
 Chicago, IL, USA, July 2015, pp.~4024--4029.

\bibitem{LiWe:15}
 D.~Lifshitz and G.~Weiss, \m Optimal control of a capacitor-type
 energy storage system, {\em IEEE Trans. on Automatic Control},
 vol.~60, 2015, pp.~216-220.

\bibitem{NaWe:15}
 V.~Natarajan and G.~Weiss, \m Almost global asymptotic stability
 of a grid-connected synchronous generator, \m submitted in 2015.

\bibitem{NaWe:16}
 V.~Natarajan and G.~Weiss, \m Synchronverters with better stability
 due to virtual inductors, virtual capacitors and anti-windup, \m
 submitted in 2016.

\bibitem{Orte:1998}
 R.~Ortega, A.~Lor\'ia, P.J.~Nicklasson and H.~Sira-Ram\'irez, \m
 {\em Passivity-Based Control of Euler-Lagrange Systems}, 
 Springer-Verlag, London, 1998.

\bibitem{OrRo:12}
 R.~Ortega and J.G.~Romero, \m Robust integral control of
 port-Hamiltonian systems: The case of non-passive outputs with
 unmatched disturbances, {\em Systems \& Control Letters}, vol.~61,
 2012, pp.~11-17.

\bibitem{Sont:1990}
 E.D.~Sontag, {\em Mathematical Control Theory: Deterministic 
 Finite Dimensional Systems}, \m Springer-Verlag, New York, 1990.

\bibitem{WeNa:16}
 G.~Weiss and V.~Natarajan, \m Integral control of stable nonlinear 
 systems, \m submitted in 2016.

\bibitem{WuLu:04}
 F.~Wu and B.~Lu, \m Anti-windup control design for exponentially
 unstable LTI systems with actuator saturation, {\em Systems \& 
 Control Letters}, vol.~52, 2004, pp.~305-322.

\bibitem{ZaccTeel:02}
 L.~Zaccarian and A.R.~Teel, \m A common framework for anti-windup,
 bumpless transfer and reliable designs, {\em Automatica}, vol.~38,
 2002, pp.~1735-1744.
}
\end{thebibliography}
\end{document}